\documentclass[12pt,reqno]{amsart}

\usepackage{amsmath,amsthm,amscd,amsfonts,amssymb,graphicx,color}

\usepackage[bookmarksnumbered,colorlinks,plainpages]{hyperref}
\hypersetup{colorlinks=true,linkcolor=red,anchorcolor=green,citecolor=cyan,urlcolor=red,filecolor=magenta,pdftoolbar=true}

\usepackage{mathrsfs}

\newtheorem{theorem}{Theorem}[section]
\newtheorem{lemma}[theorem]{Lemma}
\newtheorem{proposition}[theorem]{Proposition}
\newtheorem{corollary}[theorem]{Corollary}
\theoremstyle{definition}
\newtheorem{definition}[theorem]{Definition}
\theoremstyle{remark}
\newtheorem{remark}[theorem]{Remark}

\numberwithin{equation}{section}

\begin{document}

\setcounter{page}{1}

\title[A boundary formula]{
A boundary formula for reproducing kernel \\ Hilbert spaces of  real harmonic\\  functions  in Lipschitz domains}

% -----

\author[A. Chaira and S. Touhami ]{
 Abdellatif Chaira $^1$ and Soumia Touhami $^2$}

\address{Universit\'{e} Moulay Ismail, Facult\'{e} des Sciences, 
Laboratoire de Math\'{e}matiques et leures Applications, 
Equipe EDP et Calcul Scientifique, BP 11201 Zitoune, 50070 Mekn\`{e}s, Morocco.}

\email{\textcolor[rgb]{0.00,0.00,0.84}{a.chaira@fs.umi.ac.ma $^{1}$; touhami16soumia@gmail.com $^{2}$}}

% -----

% -----------------------------------

\keywords{ reproducing kernel Hilbert spaces, Lipschitz domains, harmonic spaces, trace spaces  and Moore--Penrose pseudo-inverse.}

% -----------------------------------

\begin{abstract}
This paper develops a new Hilbert space method to characterize a family of reproducing kernel Hilbert spaces of real harmonic functions in a bounded Lipschitz domain  $\Omega \subset \mathbb R^d, d\geq 2$  involving some families of positive self-adjoint operators and making use of  characterizations of their trace data  and of a special inner product on $H^1(\Omega).$ We also establish  boundary representation results for this family in terms of  the $L^2-$ Bergman kernel.  In particular, a boundary integral representation for the very weak solution of the  Dirichlet problem for Laplace's equation with $L^2- $ boundary data  is provided. Reproducing kernels and  orthonormal bases for the harmonic spaces are  also found. 
\end{abstract}

\maketitle

% -------------------------------------------

\section{Introduction}
We mean by a reproducing kernel Hilbert space (RKHS)  a Hilbert space associated with a kernel such that every evaluation functional is bounded. 
The concept of reproducing kernel was used for the first time on boundary value problems for harmonic and biharmonic functions by Zaremba in 1907 \cite{Z, Z2},  where he introduced the kernel corresponding to a class of functions which satisfy reproducing property and his idea to apply the kernels to the solution of boundary value problems was developed by Bergman and Schiffer \cite{BS} where the kernels were proved powerful tool for solving boundary value problems of partial differential equations of elliptic type. 
In 1950, this notion  was studied  in a systematic way by Aronszajn  \cite{Ar,Ar1} and Bergman, and was developed also by Schwartz in \cite{Sch}. 
Some results on reproducing kernels for solutions of second-order elliptic boundary
value problems were  described by Bergman and Schiffer \cite{BS} and 
significant developments of the theory of reproducing kernel Hilbert spaces were established. Several studies on reproducing kernels for real Hilbert spaces of solutions of partial differential equations do exist as well. J-L. Lions in \cite{L1}and \cite{L2} and Englis, Lukkausen, Peetre and Persson in \cite{ELPP}  described characterizations of the reproducing kernel on various different spaces of real harmonic functions where they required considerable smoothness of the boundary and involved the use of classical Green’s functions and eigenfunctions of the Laplace-Beltrami operator on the boundary of the domain. The Laplace's equation was a prototype in their study and they  extended their results to large classes of elliptic equations. In \cite{Au1}, Auchmuty followed a different approach to characterize a family of real harmonic functions Hilbert spaces  in  a bounded domain  $\Omega \subset \mathbb R^d$ and whose boundary values belong to $H^s(\partial \Omega)$ for $s\geq 0$ where  he required the boundary satisfies Gauss-Green, Rellich and compact trace theorems.  In his paper, he discussed the values of $s$ for which  they are reproducing kernel Hilbert spaces and provided an explicit  formula for the associated  reproducing kernels. His approach depended on results about the harmonic Steklov eigenfunctions of the domain. 
In \cite{Au2}, Auchmuty described the singular value decomposition of the Poisson kernel for the Dirichlet problem for the Laplace's equation in a bounded domain of $\mathbb R^d, d\geq 2.$ Under the same assumptions on the domain as in \cite {Au1}, he   obtained  a continuity result for the normal derivative operator. His paper contains  a characterization of the $L^2-$ Bergman space and a description for its associated reproducing kernel. The present work deals with a different approach to describe representation results for the family of reproducing kernel Hilbert spaces   $\mathcal H^s(\Omega)= \{ v\in H^s(\Omega) \ | \  \Delta v=0  \},$ of real harmonic functions (in the weak sense) on the usual Sobolev space $H^s(\Omega)$   
 in a bounded Lipschitz domain $\Omega \subset \mathbb R^d, d\geq 2$   for the range of values $0 \leq s <\frac{3}{2},$ where  $\Delta$ is the Laplacian. This approach also makes it possible to derive a boundary formula for the very weak solution   of the Dirichlet problem for the Laplace's equation with $L^2-$ boundary data. Other related results are also established.
% -------------------------------------------
\section{main results}
\label{results}
Let $\Omega$ be a bounded Lipschitz domain in $\mathbb R^d, d\geq 2$ with boundary $\partial \Omega$ and closure $\overline{\Omega}.$ To describe our main results we need to fix some notation. Let  $$ \mathcal H^s(\Omega) = \{  v\in H^s(\Omega) \ / \ \Delta v =0 \ in \ \mathscr{D}' (\Omega) \} $$ be the space of real harmonic functions on the usual Sobolev space $H^s(\Omega),$  where $\mathscr{D}' (\Omega)$ denotes the distributions space on $\Omega,$ and let $\Gamma$
be the trace operator from $H^1(\Omega)$ to $L^2(\partial \Omega)$ and $\Lambda$
 its Moore-Penrose inverse, that is, for $v \in H^1(\Omega),$ $\Lambda \Gamma v$ is the harmonic extension to $\Omega$   of $\Gamma v,$ where $H^1(\Omega)$ is induced by the following inner product 
  \begin{equation}
  \label{SIP}
  (u,v)_{\partial, \Omega}= \int_{\Omega} \nabla u \nabla v \ dx 
  + \int_{\partial \Omega} \Gamma u \Gamma v \ d \sigma 
  \quad \forall u,v  \in H^1(\Omega),
  \end{equation} 
   where  $\nabla v$ is the gradient of the function $v$   and $\sigma$ is the surface measure  on the boundary $\partial \Omega.$
  We denote by $ \mathcal H_{s}(\Omega)$ the following: $$ \mathcal H_{s}(\Omega) = \{ (I+\Lambda \Lambda^*)^{-(s-1)} v \ / \ v \in \mathcal H^1(\Omega) \}  \ \ \ \mbox{for} \ s\geq 1,$$ where $\Lambda^*$ denotes the adjoint operator of $\Lambda.$
Our first main result in this paper provides a functional characterization for $\mathcal H^s(\Omega)$ when $1 \leq s <3/2.$
 \begin{theorem}
 \label{thm:first main result}
   Let $\Omega$ be a  bounded Lipschitz domain of $\mathbb R^d,d \geq 2.$ Then, for all $1 \leq s <3/2,$ we have
    $$ \mathcal H^s(\Omega)= \mathcal H_{s}(\Omega) $$ with equivalence of norms. 
  \end{theorem}Consider now the embedding operator $E$ of $H^1(\Omega)$ in $L^2( \Omega).$ Its adjoint operator  denoted  $E^*,$ with respect to the inner product \eqref{SIP} on $H^1(\Omega)$, is  the solution operator of the following Robin problem for the Poisson equation
 \begin{equation}
  \label{eq:1}
          \begin{cases}
            -\Delta u=f & \text{  }  (\Omega) \\
            \partial_{\nu}u+\Gamma u =0 & \text{}  (\partial \Omega),
          \end{cases}
          \end{equation}
          where $f\in L^2(\Omega)$ and $\partial_{\nu}$ is the normal derivative operator with exterior normal $\nu.$
Denote by $E_0^*$  the solution operator of the following  Dirichlet problem for the Poisson equation
 \begin{equation}
 \label{eq:2}
              \begin{cases}
                -\Delta u^0=f & \text{  }  (\Omega) \\
                \Gamma u^0=0 & \text{}  (\partial \Omega).
              \end{cases}
              \end{equation}
 By setting $E_1^*= E^*-E_0^*$ and $u^1 =E_1^*f,$ it follows that $u^1$ is the solution of the following  Dirichlet problem for the Laplace's equation
  \begin{equation}
  \label{eq:DR}
               \begin{cases}
                 \Delta u^1=0 & \text{  }  (\Omega) \\
                \Gamma u^1=\Gamma u & \text{}  (\partial \Omega).
               \end{cases}
               \end{equation} 
                The use of the following Rellich--Ne\v{c}as lemma  allows to define a very weak solution to the following Dirichlet problem for Laplace's equation        
                  \begin{equation}
                      \label{eq:7}
                              \begin{cases}
                                \Delta v=0 & \text{  }  (\Omega) \\
                                v =g & \text{}  (\partial \Omega),
                              \end{cases}
                              \end{equation}
                     where $g\in L^2(\partial \Omega).$ The concept of a very weak solution of (\ref{eq:7}) goes back to Ne\v cas \cite[chapter 5]{N}, it is a special kind of distributional solution.
  \begin{lemma}{\cite[Chapter 5]{N}}
          Let $f\in L^2(\Omega)$ and let $u^0=E_0^*f$ be the solution of the Dirichlet problem 
          for the Poisson equation \eqref{eq:2}. Then, $\partial _{\nu} u^0 \in L^2(\partial \Omega)$. 
          Moreover, there exists a constant $c_{\Omega}>0$ depending on the geometry of $\Omega$, such that
          $$  
          \|\partial_{\nu} u^0 \|_{0,\partial \Omega} \leq  c_{\Omega} \ \|f\|_{0,\Omega},
          $$
          where $ \|\cdot \|_{0, \Omega}$ and  $ \|\cdot \|_{0,\partial \Omega}$ denote the norms on $L^2(\Omega)$ and $L^2(\partial \Omega)$ respectively.
          \end{lemma}    
           Indeed, we say that $v\in L^2(\Omega)$ is a \textbf{very weak solution} of the problem \eqref{eq:7} if for all $u$ in $H_{\Delta}^1(\Omega) \cap H_0^1(\Omega),$ we have 
            \begin{equation*}
            - \int_{\Omega} v \ \Delta u \ dx + \int_{\partial \Omega} g \  \partial _{\nu} u \ d\sigma  =0, 
            \end{equation*}
   where    $
   H_{\Delta}^1(\Omega)
   = \left\{  v \in H^1(\Omega) \ / \ \Delta v \in L^2(\Omega) \right\}
   $ and $H_0^1(\Omega)$ is the closure in 
      $H^1(\Omega)$ of infinitely differentiable functions 
      compactly supported in $\Omega$ and where the normal derivative operator $\partial_{\nu}$ is considered from $H_{\Delta}^1(\Omega)$ to $L^2(\partial \Omega).$
 Denote by $K$ the solution operator of the problem \eqref{eq:7} and by $K^*$ its adjoint operator. The following regularity  result provides a decomposition for the solution operator of the Dirichlet problem for the Laplace's equation \eqref{eq:DR}. This result will be particularly useful to characterize $\mathcal H^s(\Omega)$ for $0 \leq s\leq 1.$
      \begin{theorem} 
      \label{thm: regularity for K}
      Let $E_1^*$ be defined as above, let $\Gamma^*$ the adjoint operator of $\Gamma,$ and let   $\mathcal R(E_1^*)$ and  $\mathcal R(\Gamma ^*)$ be their range spaces  respectively. Then, $\mathcal R(E_1^*) \subset \mathcal R(\Gamma ^*).$ Moreover,   the following decomposition $  E_1^* = \Gamma ^* K^*$  holds.
      \end{theorem}
   Let us now set $\Gamma_0^* = F_1^* (I+F_1 F_1^*)^{-1/2} \Gamma ^* ,$
      where $F_1^*$ is the Moore-Penrose inverse of $E_1^*$ with adjoint operator $F_1. $  In the following theorem, an orthonormal basis for  the  $L^2-$Bergman space denoted $\mathcal H(\Omega)$ is found.
       \begin{theorem} 
               \label{thm:comp}
              Let $K$ be the solution operator of the Dirichlet problem for the Laplace's equation \eqref{eq:7} and let $K^*$ be its adjoint operator. Then, we have
                $$
                \Gamma_0^*  K^*=(I+F_1 F_1^*)^{-1/2} P_{\mathcal H(\Omega)}, $$ where $P_{\mathcal H(\Omega)}$ is the orthogonal projection onto  $\mathcal H(\Omega).$ Moreover,  $\Gamma_0^* K^*$  is a compact self-adjoint operator and there exists a sequence of  pairs $ ((\kappa_n, \phi_n)) _{n\geq 1}$ in $\mathbb R _+ ^*\times \mathcal H(\Omega) $ such that $\kappa_n$ goes to zero when $n$ goes to $+\infty$  and that for all $n \geq 1,$ one has
                       $$ \Gamma _0^* K^* \phi _n =  \kappa _n^2 \phi_n $$
               and $(\phi_n)_{n\geq 1} $ is an orthonormal basis for $\mathcal H(\Omega).$
              \end{theorem}
    Let us consider for $s\geq 0$ the family 
          $\mathcal X^s(\Omega) = \{ (I+F_1^*F_1)^{-s/2} v  \ | \ v \in \mathcal H(\Omega) \}.$  Functional and spectral characterizations for $\mathcal H^s(\Omega)$ when $ 0\leq s \leq 1$ are given in the following theorem.
         \begin{theorem}
         \label{thm: main result}
             Let $\Omega\subset \mathbb R^d, d\geq 2$ be a bounded Lipschitz domain. Then for $0\leq s\leq 1,$
          $\mathcal H^s(\Omega) = \mathcal X^s(\Omega)$ with equivalence of norms. Moreover, for all $v\in \mathcal H^s(\Omega),$
              $\displaystyle{\sum _{n}} \frac{1}{\kappa_n ^{4s}} \ |(v,\phi_n)_{0,\Omega}|^2 $ converges, where  $((\kappa_n, \phi_n))_{n\geq 1} $ is the sequence of pairs in $\mathbb R_+ ^* \times \mathcal H(\Omega)$ associated with $\Gamma_0^* K^*$ stated in Theorem \ref{thm:comp}. 
            \end{theorem}
       Let us consider now  the orthonormal basis  $(\phi_n)_{n\geq 1} $  for $\mathcal H(\Omega)$ given in Theorem \ref{thm:comp}. The Bergman kernel associated with the $L^2-$Bergman space $\mathcal H(\Omega)$ is given by  
       \begin{equation}
       \label{bergman kernel}
        b(x,y) = \displaystyle \sum_{n=1}^{\infty} \phi_n(x) \phi_n(y).
       \end{equation}
       The following theorem provides a boundary formula for an arbitrary bounded operator from $L^2(\partial \Omega)$ to $  L^2( \Omega) $ with range $\mathcal R(A)$ such that $\mathcal R(A) \subset \mathcal H(\Omega).$ 
       \begin{theorem}
       \label{thm:BF}
               Let $\Omega \subset \mathbb R^d, d\geq 2$ be a bounded Lipschitz domain, let $A$ be a bounded operator from $L^2(\partial \Omega)$ to $  L^2( \Omega) $  such that  $\mathcal R(A) \subset \mathcal H(\Omega),$ and let $A^*$ be its adjoint operator. Then, for all $v\in \mathcal R(A)$ there exists $g\in L^2(\partial \Omega)$ such that
           $$v(x)=\int_{\partial \Omega} \Big(  A^* b_x \Big ) (y) g(y) d\sigma(y),$$
           where $b(\cdot,\cdot)$ is the Bergman kernel given by \eqref{bergman kernel}.
                     \end{theorem}

   This leads to a particular boundary formula for $\mathcal H^s(\Omega)$ for the range of values $1/2 <s<3/2$ as follows:
    \begin{corollary}
    \label{co:1}
      Assume $1/2 <s < 3/2$ and let  $b(\cdot,\cdot)$ be the Bergman kernel given by \eqref{bergman kernel}. Then for all $v\in \mathcal H^s(\Omega),$ its embedding $\tilde{v} $ in $L^2(\Omega)$ has the following form
       
    $$
          \tilde{v}(x)= \int_{\partial \Omega} \Big(  K^* b_x \Big ) (y) \Big( \widehat{\Gamma}_s v \Big )(y) \ d\sigma(y),
         $$ 
       where $\Gamma_s$ denotes the trace operator from $H^s(\Omega)$ to $H^{s-1/2}(\partial \Omega,)$   and $\widehat{\Gamma}_s v $ is the embedding of  $\Gamma_s v$   in $L^2(\partial \Omega).$ 
    \end{corollary}

In particular, a boundary representation result for the very weak solution of the Dirichlet problem for the Laplace's equation \eqref{eq:7} is derived:
\begin{corollary}
\label{co:2} Let $K$ be the solution operator of the problem \eqref{eq:7} with range $\mathcal R(K).$ Then,
 for all $v\in \mathcal R(K),$ there exists $g\in L^2(\partial \Omega)$ such that
   
    \begin{equation}
    \label{Lions's formula}
    v(x)=\int_{\partial \Omega} \Big(  K^* b_x \Big ) (y) g(y) d\sigma(y),
    \end{equation}
    where  $b(\cdot,\cdot)$ is the Bergman kernel given by \eqref{bergman kernel}.
      \end{corollary}
      A boundary formula for the operator $\Gamma_0^*$ is provided as well:
      \begin{corollary}
      \label{co:3} Let $\Gamma_0^*$ defined as above with range $\mathcal R(\Gamma_0^*)$ and let $\Gamma_0$ be its adjoint operator.  Then,
       for all $v\in \mathcal R(\Gamma_0^*),$ there exists $z\in L^2(\partial \Omega)$ such that
         
          \begin{equation}
          v(x)=\int_{\partial \Omega} \Big(  \Gamma_0 b_x \Big ) (y) z(y) d\sigma(y),
          \end{equation}
          where  $b(\cdot,\cdot)$ is the Bergman kernel given by \eqref{bergman kernel}.
            \end{corollary}
  An orthonormal basis for $\mathcal H^s(\Omega), 0\leq s \leq 1$ is constructed and a description of its reproducing kernel is given.
   \begin{theorem}
    \label{prop:repro}
   Let $((\kappa_n, \phi_n))_{n\geq 1} $ in $\mathbb R_+ ^* \times \mathcal H(\Omega)$ the  sequence of  pairs associated with the operator $\Gamma_0^* K^*$ considered in Theorem \ref{thm:comp}. Then, for $0 \leq s \leq 1$ $((\kappa_n^{2s}\phi_n))_{n\geq 1}$ is an orthonormal basis for  $\mathcal H^s(\Omega).$ Moreover,
                 the associated  reproducing kernel denoted  $b^s(\cdot,\cdot)$  is given by:
                   $$b^s(x,y) = \displaystyle \sum_{n=1}^{\infty} \phi_n^s(x) \phi_n^s(y) $$
       where $\phi_n^s = \kappa_n^{2s} \phi_n .$
      \end{theorem}  
 This paper is organized as follows. In section ~\ref{sec:02} we present some general tools and basic results  which will be  frequently used in the sequel and we recall  some results on the Moore-Penrose inverse,  also the description of various function spaces including  Sobolev spaces in Lipschitz domains and the main ingredients needed  to characterize the space of real harmonic functions on $H^s(\Omega)$ for $0\leq s < 3/2.$ We will make frequent use of the results in our previous paper \cite{TCT} where a new  contribution concerning the Moore-Penrose inverse is given  and  functional characterizations for the trace spaces $H^s(\partial \Omega)$ for the range of values $-1\leq s \leq 1$ are provided. The main developments in this paper and proofs of the main results are encountered in section \ref{section:main one}.

\section{Background material}
\label{sec:02}

Let $\mathcal H_1$ and $\mathcal H_2$ be two Hilbert spaces with inner products 
$(\cdot,\cdot)_{\mathcal H_1}$ and $(\cdot,\cdot)_{\mathcal H_2}$ and associated 
norms $\|\cdot\| _{\mathcal H_1}$ and $\|\cdot\| _{\mathcal H_2}$, respectively.  
We need  first to fix some notations. By $\mathcal L(\mathcal H_1, \mathcal H_2)$, we denote 
 the space of all linear operators from $\mathcal H_1$ 
into $\mathcal H_2$ and $\mathcal L(\mathcal H_1, \mathcal H_1)$ 
is briefly denoted by $\mathcal L(\mathcal H_1)$. For an operator 
$A\in \mathcal L(\mathcal H_1,\mathcal H_2)$, 
$\mathcal{D}(A)$, $\mathcal R(A)$ and $\mathcal N(A)$ 
denote its domain, its range and its null space, respectively. 
For $A, B \in \mathcal L(\mathcal H_1, \mathcal H_2)$,  
$B$ is called an extension of $A$ if $ \mathcal  D(A)\subset  \mathcal D(B)$ 
and $Ax  =  Bx$ for all $x\in \mathcal D(A)$, and this fact is denoted by  
$A\subset B$. The set of all bounded operators from $\mathcal H_1$ 
into $\mathcal H_2$ is denoted by $\mathcal B(\mathcal H_1, \mathcal H_2)$, 
while $\mathcal B(\mathcal H_1, \mathcal H_1)$ is briefly denoted 
by $\mathcal B(\mathcal H_1)$. The set of all closed densely defined 
operators from $\mathcal H_1$ into $\mathcal H_2$ is denoted 
by $\mathcal C(\mathcal H_1, \mathcal H_2)$, and 
$\mathcal C(\mathcal H_1, \mathcal H_1)$ is denoted by
$\mathcal C(\mathcal H_1)$. For  $A \in \mathcal C(\mathcal H_1,\mathcal H_2)$, 
its adjoint operator is denoted by $A^*\in \mathcal C(\mathcal H_2,\mathcal H_1)$.  
A self-adjoint operator $A$ on a Hilbert space $\mathcal H$ 
is said to be positive (resp., strictly positive) if $(Ax,x)_{\mathcal H} \geq 0$ for all  
$x \in \mathcal D(A)$ ($(Ax,x)_{\mathcal H} > 0$); in such case we write $A\geq 0$ (resp., $A> 0$).
The next theorem goes back to Douglas \cite{Do} and it will prove to be useful in this paper.
\begin{theorem}{\cite[Theorem 1]{Do}} 
\label{thm:DT}
Let $\mathcal H$ be a Hilbert space, and let $A,B \in \mathcal B(\mathcal H)$ 
be two bounded operators. The following statements are equivalent:
\begin{enumerate}
\item $\mathcal R(A)\subset \mathcal R(B)$,
\item $AA^* \leq \mu \ BB^*$ for some $\mu \geq 0$,
\item there exists a bounded operator $C\in \mathcal B(\mathcal H)$ such that $A=BC$.
\end{enumerate}
Moreover, if the above items 1, 2 and 3 hold, then there exists 
a unique operator $C\in \mathcal B(\mathcal H)$ such that
\begin{enumerate}
\item[(a)] $\|C\|^2= \inf \{ \mu \ | \ AA^* \leq \mu \  BB^* \}$,
\item[(b)] $\mathcal N(C) =\mathcal N(A)$,
\item[(c)] $\mathcal R(C) \subset \overline{\mathcal R(B^*)}$,
\end{enumerate}
where $\overline{\mathcal R(B^*)}$ denotes the closure of $\mathcal R(B^*).$
\end{theorem}

\begin{definition}
Let ${\mathcal H}_1$ and ${\mathcal H}_2$ be two Hilbert spaces, 
$A \in {\mathcal C}({\mathcal H}_1, {\mathcal H}_2)$ be a closed densely 
defined operator and $A^*$  its adjoint. The Moore--Penrose 
inverse of $A$, denoted by $A^\dagger$, is defined as the unique 
linear operator in  ${\mathcal C }({\mathcal H}_2, {\mathcal H}_1)$  
such that
$$
{\mathcal D}(A^\dagger)={\mathcal R}(A)
\oplus {\mathcal N}(A^*), 
\quad  \mathcal N(A^{\dagger })= \mathcal N(A^*) 
$$ 
and 
\begin{equation*}
\begin{cases}
AA^\dagger A=A, &\text{  }  \\
A^\dagger AA^\dagger=A^\dagger, &\text{}
\end{cases}
\qquad
\begin{cases}
AA^\dagger \subset P_{\overline{{\mathcal R}(A)}}, & \text{  }  \\
A^\dagger A \subset P_{\overline{{\mathcal R}( A^\dagger)}}, &\text{}
\end{cases}
\end{equation*}
where $P_{{\mathcal E}}$ denotes the orthogonal projection 
on the closed subspace ${\mathcal E}.$ 
   \end{definition}

A  result due to Von Neumann (see \cite{Gro,L}) asserts that 
for $A\in \mathcal C(\mathcal H_1,\mathcal H_2)$ the operators 
$(I+AA^*)^{-1}$ and $A^*(I+AA^*)^{-1}$ are everywhere defined and bounded and that $(I+AA^*)^{-1}$ is self-adjoint.
Similarly, the operators $(I+A^*A)^{-1}$ and $A(I+A^*A)^{-1} $ are everywhere 
defined and bounded, and $(I+A^*A)^{-1}$ is self-adjoint. 
Moreover, 
\begin{equation*}
(I+AA^*)^{-1} A \subset A(I+A^*A)^{-1}
\end{equation*}
and
\begin{equation*}
(I+A^*A)^{-1} A^*  \subset A^*(I+AA^*)^{-1}.
\end{equation*}

\begin{proposition}[See Lemma 2.5 and Corollary 2.6 of \cite{L}]  
\label{prop:3.1}
Let $A \in \mathcal C(\mathcal H_1, \mathcal H_2)$ 
and $B \in \mathcal C(\mathcal H_2,\mathcal H_1)$ be
such that $ B=A^{\dagger}$. Then,
\begin{enumerate}
\item $A(I+A^*A)^{-1} = B^*(I+BB^*)^{-1}$,
\item $(I+A^*A)^{-1}+(I+BB^*)^{-1}= I+P_{\mathcal N(B^*)}$,   
\item $A^*(I+AA^*)^{-1} = B(I+B^*B)^{-1}$,
\item $(I+AA^*)^{-1}+(I+B^*B)^{-1}= I+P_{\mathcal N(A^*)}$,
\item $(I+AA^*)^{-1}+(I+B^*B)^{-1}= I$  (if $A^*$ is injective),
\item $\mathcal N(A^*(I+AA^*)^{-1/2}) =\mathcal N(A^*) = \mathcal N(B)$.
\end{enumerate}
\end{proposition}

\begin{lemma}[See Proposition~1.7 of \cite{LM}] 
\label{lem:LM}
Let  ${\mathcal H}_1, {\mathcal H}_2$ be two Hilbert spaces, let
$A \in {\mathcal B}({\mathcal H}_1, {\mathcal H}_2),$ 
and let $B$ be its Moore--Penrose inverse. Then,
one has 
\begin{equation*}
\|x\|_{\mathcal H_1}^2= \|B^*(I+BB^*)^{-1/2} 
x\|_{\mathcal H_2}^2 + \|(I+BB^*)^{-1/2} x\|_{\mathcal H_1}^2
\end{equation*}
for all $x\in \mathcal H_1$.
\end{lemma}

\begin{theorem}   [See Theorem~3.5 of \cite{TCT}] 
\label{thm:3.5}
Let  ${\mathcal H}_1$ and ${\mathcal H}_2$ be two Hilbert spaces, let $A \in {\mathcal B}({\mathcal H}_1, {\mathcal H}_2)$ and let
$B$ be its Moore--Penrose inverse. Then, the operator $B^*(I+BB^*)^{-1/2}$ 
is bounded with closed range and has a bounded Moore--Penrose inverse given by 
\begin{equation*}
T_B=B(I+B^*B)^{-1/2}+A^*(I+B^*B)^{-1/2}.
\end{equation*}
Moreover, the adjoint operator of $T_B$ is $T_{B^*}$, where 
\begin{equation*}
T_{B^*}=B^*(I+BB^*)^{-1/2}+A(I+BB^*)^{-1/2}.
\end{equation*}
\end{theorem}

\begin{corollary} [See Corollary~3.6 of \cite{TCT}] 
\label{cor:3.6}
Let $A\in\mathcal B(\mathcal H_1,\mathcal H_2)$ 
and let $B$ be its Moore--Penrose inverse. Then,
$B^*(I+BB^*)^{-1/2}$ is an isomorphism of
$\mathcal N(B^*)^{\perp}$ and $\mathcal R(B^*),$ where $\mathcal N(B^*)^{\perp}$  is the orthogonal complement of $\mathcal N(B^*).$ 
\end{corollary}
\begin{corollary} [See Corollary~3.7 of \cite{TCT}] 
\label{cor:2.8}
Let $A\in\mathcal B(\mathcal H_1,\mathcal H_2),$ and let $B$ be 
its Moore--Penrose inverse. Then, $T_B$ is an isomorphism of
$\mathcal R(B^*)$ and $\mathcal N(B^*)^{\perp}$.  
\end{corollary}

\begin{theorem}[See Theorem 3.8 of \cite{TCT}] 
\label{thm:3.8}
Let $\mathcal H_1$ and $\mathcal H_2$ be two Hilbert spaces, let
$A\in \mathcal B(\mathcal H_1,\mathcal H_2),$ and let $B$ be 
its Moore--Penrose inverse. Then, the decomposition
$$
A= (I+B^*B)^{-1/2} T_{B^*}
$$
holds, where $T_{B^*}=B^*(I+BB^*)^{-1/2}+A(I+BB^*)^{-1/2}$.
\end{theorem}
      
Another theoretical background we should include is Reproducing Kernel Hilbert Spaces.
   A Hilbert space $\mathcal H$ of functions defined on an open subset $\Omega \subset \mathbb{R}^d, d\geq 1$ is said to be a reproducing kernel Hilbert space (RKHS) if $$ \forall x \in \Omega, \exists \ c _x \in \mathbb R, |f(x)|\leq c_x \|f\|_{\mathcal H} \ \ \forall f \in \mathcal H.$$
Equivalently, the evaluation functional at $x,$ $\ell_x:\mathcal H \rightarrow \mathbb R,$ \ $f \mapsto f(x)$ is a bounded linear functional.
     By Riesz representation theorem, it follows that there exists a unique function $k_x \in \mathcal H$ such that
     $$ \ell_x f = (k_x,f)_{\mathcal H} =f(x) \ \ \forall f \in \mathcal H .$$
     The function $k(.,.)$ defined by $k(x,x')=k_x(x')$ is called the \textbf{reproducing kernel} of $\mathcal H.$ Recall that a symmetric function $k(.,.): \Omega \times \Omega \longrightarrow \mathbb R$ is said to be positive definite if  $$ \forall l \in \mathbb N^*, \forall x^1,..., x^l \in \Omega, \ \forall \lambda_1,...,\lambda_l \in \mathbb R  \  \ \sum_{i,j=1}^{l} \lambda_i \lambda_j k(x^i, x^j) \geq 0.$$
It turns out that a reproducing kernel is a positive definite function. 
   A further important property is that if it exists, a reproducing kernel is unique. Equivalently, a reproducing kernel Hilbert space uniquely determines its reproducing kernel, this is the topic of the following theorem.	
    \begin{theorem}{Moore-Aronszajn \cite{Ar,Ar1}}.   
       Suppose $k$ is a symmetric, positive definite function on $\Omega,$ then there exists a unique Hilbert space of functions on $\Omega$ for which $k$ is a reproducing kernel. Conversely, given a reproducing kernel Hilbert space $\mathcal H$ then there exists a unique reproducing kernel corresponding to $\mathcal H.$
\end{theorem}
  It is also interesting to mention that if $\mathcal H$ is a reproducing kernel Hilbert space and that $(e_n)_{n\geq 1}$ is an orthonormal basis for $\mathcal H$, then the associated  reproducing kernel should take the form 
     \begin{equation}
  \label{eq: forme som}
    k(x,y)= \displaystyle{\sum_{n=1}^{+\infty}} e_n(x) e_n(y),  
   \end{equation}
(See \cite{Ar,Ar1}).
 \begin{definition}
 \label{def:LCB}
 Let $\Omega$ be an open subset of $\mathbb R^d$ with boundary 
 $\partial \Omega$ and closure $\overline{\Omega}$. 
 We say that $\partial \Omega$ is Lipschitz continuous 
 if for every $x\in \partial \Omega$ there exists a coordinate system 
 $(\widehat{y}, y_d)\in \mathbb R^{d-1}\times\mathbb R$, 
 a neighborhood $Q_{\delta,\delta'}(x)$ of $x $ and a Lipschitz function 
 $\gamma_x:\widehat{Q}_{\delta} \rightarrow \mathbb R$  with 
 the following properties:
 \begin{enumerate}
 \item $\Omega \cap Q_{\delta,\delta'}(x) 
 = \left\{(\widehat{y},y_d) \in Q_{\delta,\delta'}(x) 
 \ / \ \gamma_x(\widehat{x}) < y_{d} \right\}$;
 
 \item $\partial \Omega \cap  Q_{\delta,\delta'}(x) 
 =\left\{  (\widehat{y},y_d) \in Q_{\delta,\delta'}(x) 
 \ / \ \gamma_x(\widehat{x}) = y_{d} \right\}$;
 \end{enumerate}
 where 
 $$ 
 Q_{\delta,\delta'}(x) 
 = \left\{  (\widehat{y},y_d) \in \mathbb R^d \ / \  
 \|\widehat{y}-\widehat{x}\| _{\mathbb R^{d-1}} < \delta  
 \ \ \text{and} \ \ |y_d - x_d | < \delta' \right\}
 $$
 and 
 $$ 
 \widehat{Q}_{\delta}(x) = \{ \widehat{y} \in \mathbb R^{d-1} 
 \ / \  \|\widehat{y}-\widehat{x}\| _{\mathbb R^{d-1}} < \delta  \} 
 $$
 for  $\delta, \delta' > 0$. 
 An open connected subset $\Omega \subset \mathbb R^d$, whose boundary 
 is Lipschitz continuous, is called a Lipschitz domain.
 \end{definition}
Any Lipschitz domain  has a surface measure $\sigma$ 
     and an outward unit normal $\nu$ that exists $\sigma$-almost everywhere 
     on the boundary $\partial \Omega.$ In the rest of this article, $\Omega$ denotes a bounded Lipschitz 
 domain in $\mathbb{R}^d$, $d\geq 2$.  We denote by $\mathcal C^k(\Omega)$, $k\in \mathbb{N}$ 
  or $k= \infty$, the space of real $k$ times continuously 
  differentiable functions on $\Omega$. The space $\mathcal{C}^{\infty}$
  of all real functions on $\Omega$ with a compact support 
  in $\Omega$ is denoted by $\mathcal C_c^{\infty} (\Omega)$. For a sequence  $(\varphi_n)_{n\geq 1}$ in $\mathcal C_c^{\infty}(\Omega)$ 
   and $\varphi \in \mathcal C_c^{\infty}(\Omega)$, we say that 
   $(\varphi_n)_{n\geq 1} $ converges to $\varphi$  
   if there exists a compact $Q\subset \Omega$ such that for all 
   $n\geq 1 $ $\mathrm{supp}(\varphi_n) \subset Q$ and for all multi-index 
   $\alpha \in \mathbb{N}^d$, the sequence $(\partial^{\alpha} \varphi_n)_{n\geq 1}$ 
   converges uniformly to $\partial ^{\alpha} \varphi$ where $\partial ^{\alpha}$ denotes the partial differential derivative operation of order $\alpha.$
   The space $\mathcal C_c^\infty(\Omega)$
   induced by this convergence is denoted $\mathscr{D}(\Omega)$,
   as in the theory of distributions, 
   with $\mathscr{D}^{\prime}(\Omega)$ 
   the space of distributions on $\Omega$. We denote by $\mathscr{D}(\overline{\Omega})$  the set of all restrictions on $\overline{\Omega}$ of functions in $\mathcal C_c^\infty(\mathbb R^d)$ and by $\mathcal C^1(\overline{\Omega})$ the set of all restrictions on $\overline{\Omega}$ of functions in $C^1(\mathbb R^d).$ For real $s,$ consider the usual Sobolev  spaces $H^s(\Omega)$   defined on $\Omega$ with  usual inner product and  norm  denoted $(\cdot, \cdot)_{s,\Omega}$ and $\|\cdot\|_{s,\Omega}$ respectively. The trace spaces   $H^s(\partial \Omega)$ are defined on $\partial \Omega$ for $-1 \leq s \leq 1$  with usual inner product and  norm  denoted  by  $(\cdot, \cdot)_{s,\partial \Omega}$ and $\|\cdot\|_{s,\partial \Omega}$ respectively.
 The trace operator maps each continuous function $u$ on $\overline{\Omega}$ 
 to its restriction onto $\partial\Omega$ 
 and may be extended to be a bounded surjective operator, 
 denoted by $\Gamma_s$, from $H^s(\Omega)$ to 
 $H^{s-1/2}(\partial\Omega)$ for $1/2<s<3/2$  
 \cite{Co,Mc} such that
 $\mathcal R(\Gamma_s)= H^{s-1/2}(\partial \Omega)$ 
and  $\mathcal N(\Gamma_s) = H_0^s(\Omega),$
 where $H_0^s(\Omega)$ is defined to be the closure in 
 $H^s(\Omega)$ of infinitely differentiable functions 
 compactly supported in $\Omega$. For $s>3/2$, the trace operator 
 from $H^s(\Omega)$ to $H^1(\partial \Omega)$ is bounded \cite{Co}. We now  set $\Gamma=T_1 \Gamma_1$, where $\Gamma_1$ is the trace 
 operator from $H^1(\Omega)$ to $H^{1/2}(\partial \Omega)$ and 
 $T_1$ is the embedding operator from  $H^{1/2}(\partial \Omega)$  
 into $L^2(\partial \Omega)$. According to a classical result 
 of Gagliardo (See \cite[Teorema 1.I]{Ga}), we know that
 $\mathcal R(\Gamma) = H^{1/2}(\partial \Omega)$. Since $\Gamma_1$ 
 is bounded and $T_1$ is compact (See \cite[Theorem 1.4.3.2]{Gris}), 
 the trace operator $\Gamma$ from $H^1(\Omega)$ 
 to $L^2(\partial \Omega)$ is therefore compact as a composition of a bounded operator and a compact one. 
 \begin{lemma}
 \label{lem:TA}
 Let $\Gamma$ be the trace operator 
 from $H^1(\Omega)$ into $L^2(\partial \Omega)$. 
 Then, $\Gamma^*$ is compact and injective.
 \begin{proof}
 The compactness of $\Gamma$ implies the compactness of $\Gamma^*$ according to Schauder's Theorem \cite{Co1}. Also, since $\mathcal R(\Gamma) = H^{1/2}(\partial \Omega)$ and that $H^{1/2}(\partial \Omega)$ is dense in $L^2(\partial \Omega)$ (See \cite[Theorem 1.4.2.1]{Gris}), $\Gamma^*$ is therefore injective.
 \end{proof}
 \end{lemma}   Now, let us induce $H^1(\Omega)$ with the inner product: 
 \begin{equation*}
 (u,v)_{\partial, \Omega}= \int_{\Omega} \nabla u \nabla v \ dx 
 + \int_{\partial \Omega} \Gamma u \Gamma v \ d \sigma 
 \quad \forall u,v  \in H^1(\Omega),
 \end{equation*} 
  where  $\nabla v$ is the gradient of the function $v.$     
 The associated norm $\|\cdot\|_{\partial, \Omega}$ is given by
 \begin{equation*}
 \|u\|_{\partial, \Omega}=\left(\|\nabla u\|^2_{0,\Omega}
 +\|\Gamma u\|^2_{0,\partial \Omega} \right)^{1/2},
 \end{equation*}
 and
 $H^{1}(\Omega)$ induced with the inner product $(\cdot,\cdot)_{\partial,\Omega}$,
 is denoted by $H_{\partial}^{1}(\Omega)$. According to Ne\v{c}as (\cite [Theorem 1.9, page 20]{N}), the norms $\|\cdot\|_{\partial, \Omega}$ 
 and $\|\cdot\|_{1, \Omega}$ are equivalent.
  The normal derivative map $\partial _{\nu} $
 transforms each  $v\in C^1(\overline{\Omega})$ on 
 $\partial_{\nu} v = \nu \cdot (\nabla v)_{| \partial\Omega}$ 
 onto $L^{\infty}(\partial\Omega)$  and  
 may be extended to be a bounded linear operator 
 $\widehat \partial _{\nu}$ from $ H_{\Delta}^1(\Omega)$ to
 $H^{-1/2}(\partial\Omega),$ where 
 $
 H_{\Delta}^1(\Omega)
 = \left\{  \  v \in H^1(\Omega) \ / \ \Delta v \in L^2(\Omega) \right\}$ and  $H^{-1/2}(\partial\Omega)$ is the dual space of $H^{1/2}(\partial\Omega).$  
 This is a consequence of the following lemma.
 
 \begin{lemma}[See Lemma~20.2 of \cite{T}]
 \label{lem:5.2}
 The application $w\longmapsto w \cdot \nu$  defined  
 from $\big (\mathscr D(\overline{\Omega}) \big ) ^d $ into 
 $L^{\infty} (\partial \Omega)$ is well defined and extends into 
 a linear continuous map from $H_{div}(\Omega)$ into the dual space 
 of $H^{1/2}(\partial\Omega)$, that is, $H^{-1/2}(\partial\Omega)$, where 
 $$
 H_{div}(\Omega)= 
 \left\{  V \in (L^2(\Omega))^d \ | \ \mathrm{div}\, V \in L^2(\Omega) \right\}
 $$
 with the norm 
 $$
 \|V\|_{div,\Omega}^2 = \sum_{i=1}^{d} \|v_i\|_{0,\Omega}^2 
 + \|\mathrm{div}\, V\|_{0,\Omega}^2 
 $$
 for any $V=(v_1, \ldots, v_d)\in  H_{div}(\Omega),$ where $\mathrm{div}\, V = \displaystyle{\sum_{i=1}^{d}}  \partial_i v_i.$
 Moreover, the mapping is surjective.
 \end{lemma}
 
 As a consequence, we have the following result.
 
 \begin{lemma}
 \label{prop:5.3}
 For all $u\in H_{\Delta}^1(\Omega)$ there exists  
 $\widehat \partial _{\nu} u  \in H^{-1/2}(\partial\Omega)$ 
 such that 
 $$ 
 \int_{\Omega} \nabla u \nabla v dx = - \int_{\Omega} \Delta u \ v  \ dx 
 + \left<\widehat{\partial _{\nu}} u, \Gamma_1 v\right> 
 $$
 for all $v\in H^1(\Omega)$, where  $\left<\cdot,\cdot\right>$ is the duality pairing 
 between $H^{1/2}(\partial\Omega)$ and $H^{-1/2}(\partial\Omega)$.
 The application $u \longmapsto \widehat{\partial _{\nu}} u$ is the continuous 
 extension of 
 $$
 u\longmapsto \nu \cdot (\nabla v)_{| \partial\Omega},
 $$  
 which is defined for all $u \in \mathscr{D}(\overline{\Omega}).$

 \end{lemma}
 
 \begin{proof}
 Let $u\in H_{\Delta}^1(\Omega)$. By setting $w=\nabla u$, we have 
 $$
 w\in \left( L^2(\Omega) \right) ^d \quad \mbox{and} 
 \quad \mathrm{div}\, w = \mathrm{div}\, \nabla u = \Delta u \in L^2(\Omega).
 $$ 
 According to Lemma~\ref{lem:5.2}, there exists 
 $w\cdot\nu \in H^{-1/2}(\partial\Omega)$ such that  
 $$ 
 \int_{\Omega} w \nabla v \ dx + \int_{\Omega} (\mathrm{div}\, w) v \ dx 
 = \left<w\cdot\nu, \Gamma_1 v \right>
 $$ 
 for all $v \in H^1(\Omega)$ or 
 $$ 
 \int_{\Omega} \nabla u \nabla v dx + \int_{\Omega} \Delta u \ v  \ dx 
 = \left<\widehat{\partial _{\nu}} u, \Gamma_1 v\right>,
 $$
 where $\left<\cdot,\cdot\right>$ is the duality pairing between 
 $H^{1/2}(\partial\Omega)$ and $H^{-1/2}(\partial\Omega)$.
 \end{proof}
 
As a consequence of Lemma~\ref{prop:5.3}, we can apply  Green's formula, as follows.
 	
 \begin{corollary}[Green's formula]
 \label{prop:GF}
 Let $\Omega$ be a bounded Lipschitz domain. Then, 
 $$ 
 \int_{\Omega} \nabla u \nabla v dx = - \int_{\Omega} \Delta u \ Ev  \ dx 
 + \left<\widehat{\partial _{\nu}} u, \Gamma_1 v\right>
 $$
 for all $u \in H_{\Delta}^1(\Omega)$  and $ v \in H^1(\Omega)$,
 where $E$ is the embedding operator from $H^1(\Omega)$ into $L^2(\Omega)$ and 
 $\left<\cdot,\cdot\right>$ is the duality pairing between 
 $H^{1/2}(\partial\Omega)$ and  $H^{-1/2}(\partial\Omega)$.
 \end{corollary}

  Now, consider the trace operator $\Gamma$ from $H_{\partial}^1(\Omega)$ 
to $L^2(\partial\Omega).$  For $g \in L^2(\partial \Omega)$, its adjoint operator
$\Gamma ^* \in \mathcal B(L^2(\partial \Omega), H_{\partial}^1(\Omega))$ is the solution operator 
of the following Laplace's equation with Robin boundary condition:
\begin{equation*}
\begin{cases}
\Delta z =0 & \text{  }  (\Omega) \\
\partial_{\nu}z+ \Gamma z= g & \text{}  (\partial \Omega),
\end{cases}
\end{equation*}
where $\partial_{\nu}$ is the normal derivative operator, considered as non-bounded, 
from $H_{\Delta}^1(\Omega)$ to $L^2(\partial \Omega)$ (see \cite[Theorem 5.5]{TCT} ). Consider also the embedding operator $E$ of  $H_{\partial}^1(\Omega)$ in $L^2(\Omega).$  For $f\in L^2(\Omega),$ its adjoint operator $E^*$ is the solution operator of the following
Poisson equation with Robin boundary condition:
\begin{equation*}
\begin{cases}
-\Delta u =f & \text{  }  (\Omega) \\
\partial_{\nu}u+ \Gamma u=0& \text{}  (\partial \Omega),
\end{cases}
\end{equation*}
(see \cite[Theorem 5.7]{TCT}).
The Moore-Penrose inverse of the trace operator denoted $\Lambda$ is the solution operator of the Dirichlet problem 
 for the Laplace's equation with data in 
 $\mathcal D(\Lambda)= H^{1/2}(\partial \Omega)$. Moreover,
 $$ 
 {\mathcal D}(\Lambda)= {\mathcal R}(\Gamma),  \quad         
 {\mathcal N}(\Lambda^*)= \mathcal N(\Gamma)= H_0^1(\Omega)
 $$
 and ${\mathcal R}(\Lambda)$ is characterized by
 \begin{equation*}
 \mathcal R(\Lambda)=
 \left\{   v\in H^{1}(\Omega) ~~~/~~ \Delta v =0 
 \ \ \mbox{in} \ \ \mathscr D '(\Omega) \right\},
 \end{equation*}
(see \cite[Theorem 5.6]{TCT}). Moreover, since $\Gamma$ is bounded, it follows from classical results of functional analysis (see \cite{Co1}), that $\mathcal R(\Lambda)$ is closed and that there exists $c^{\prime}>0$ such that for all $g\in \mathcal D(\Lambda) , $ we have $$c^{\prime}  \ \|g\|_{0, \partial \Omega} \leq \|\Lambda g\| _{\partial, \Omega} .$$ 
 The following result of McLean \cite{Mc} is a generalization  of Ne\v{c}as arguments (See \cite[Chapter 5]{N}) to strongly elliptic systems, and the version we will consider here is particularly adopted to the Laplacian and it will prove to be useful to prove Theorem~\ref{thm: regularity for K}.
\begin{theorem}{ \cite[Theorem 4.24]{Mc}}
\label{thm:MC}
Let $\Omega$ be a bounded Lipschitz domain of $\mathbb R^d$ 
and $u\in H _{\Delta}^1(\Omega)$.
\begin{enumerate}
\item If $\partial_{\nu} u \in L^2(\partial \Omega)$, 
then $\Gamma u \in H^1(\partial \Omega)$ and there exists 
a constant $c_{\Omega}>0$, depending on the geometry of $\Omega$, such that
$$ 
\| \Gamma u \|_ {1,\partial \Omega} \leq c _{\Omega} 
\  ( \|u\|_{1, \Omega}^2 + \|\Delta u\| _{0,\Omega}^2 
+ \| \partial  _{\nu}u \| _{0,\partial \Omega}^2)^{1/2}.
$$

\item If $\Gamma u \in H^1(\partial \Omega)$, then 
$\partial_{\nu} u \in L^2(\partial \Omega)$ and there exists 
a constant $c_{\Omega}^{\prime}>0$, depending on the geometry 
of $\Omega$, such that
$$ 
\| \partial_{\nu} u \|_ {0,\partial \Omega} \leq c_{\Omega}^{\prime} 
\  \left( \|u\|_{1, \Omega}^2 +\|\Delta u\| _{0,\Omega}^2 
+ \| \Gamma u\| _{1,\partial \Omega}^2\right)^{1/2}.
$$
\end{enumerate}
\end{theorem}
 We will also  make frequent use of the following result describing  the trace spaces $H^{s}(\partial \Omega)$ to handle the harmonic spaces $\mathcal H^{s}( \Omega)$ in Section \ref{section:main one}.
\begin{theorem}[See Corollary~6.9 of \cite{TCT}] 
\label{thm:trace}
Let $\Gamma$ be the trace operator from $H_{\partial}^1(\Omega)$ 
to $L^2(\partial \Omega)$, let $\Lambda $ be its Moore--Penrose 
inverse, and let $\Lambda ^*$ be its adjoint operator. Then, for $0\leq s \leq 1,$ $\mathcal H^{s}(\partial \Omega)= \{  (I+\Lambda ^* \Lambda)^{-s}g \ | \ g\in L^2(\partial \Omega)\}$ form an interpolating family
and $\mathcal H^{s}(\partial \Omega)
=  H^{s}(\partial \Omega)$ with equivalence of norms.

\end{theorem}

% -------------------------------------------

\section{Proofs of the main results}
\label{section:main one}
In this section, we give detailed proofs of the main results stated in Section \ref{results}.  Let $\Omega$ be a bounded Lipschitz domain of $\mathbb R^d, d\geq 2,$  let $\Gamma \in {\mathcal B}( H_{\partial}^1(\Omega),L^2(\partial \Omega))$ be the trace operator,  and let 
$\Lambda \in {\mathcal C}(L^2(\partial \Omega), H_{\partial}^1(\Omega)) $ be its Moore-Penrose inverse. For real $s\geq 0$ consider the family of Hilbert  spaces  
  $$ \mathcal H^s(\Omega) = \{  v\in H^s(\Omega) \ / \ \Delta v =0 \ in \ \mathscr{D}' (\Omega) \} $$ of real harmonic functions on the Sobolev space $H^s(\Omega).$  We will establish new characterizations of the family $\mathcal H^s(\Omega)$  for the range of values $ 0\leq s < 3/2$ involving some families of positive self-adjoint  operators. Two cases will be discussed. Also, boundary formulas for $\mathcal H^s(\Omega)$ when $1/2<s<3/2$ and for the very weak solution of the Dirichlet problem for the Laplace's equation \eqref{eq:7} will be derived.
\subsection{ \textbf {The case $ 1\leq s < 3/2$}}

Consider for $s\geq 0$ the family of Hilbert spaces $\mathcal H^s(\partial \Omega)=  \left\{ (I+\Lambda^* \Lambda)^{-s} g 
\ | \ g\in L^2(\partial \Omega) \right\}$ and for $g,y \in \mathcal H^s(\partial \Omega),$  consider  the norm 
  $$\|y\|_{s,\partial}= \|y\|_{\mathcal H^s(\partial  \Omega)} = \| (I+\Lambda ^* \Lambda)^{s}y\| _{0,\partial \Omega},$$
associated with the inner product 
  $$(g,y)_{s,\partial} =   \left((I+\Lambda ^* \Lambda) ^{s}g, (I+\Lambda ^* \Lambda)^{s}y \right)_{0,\partial \Omega} . $$

By considering the following notation
  $$ \mathcal H_{s}(\Omega) = \{ (I+\Lambda \Lambda^*)^{-(s-1)} v \ / \ v \in \mathcal H^1(\Omega)  \}  \ \ \ \mbox{for} \ s\geq 1,$$
the first task in this subsection will be to describe the relationship between $\mathcal H^s(\Omega)$ and $\mathcal H_{s}(\Omega) $ for $ 1 \leq s < 3/2.$ To this end, consider the trace operator $\Gamma_s$ from $ H^s(\Omega)$ to $ H^{s-1/2}(\partial \Omega), $  it is an isomorphism of $\mathcal H^s(\Omega)$ and $ H^{s-1/2}(\partial \Omega)$ (see \cite{Co} and \cite{Mc}), which implies the existence of  real positive constants  $c_1$ and $c_2$ such that
 for all $u\in \mathcal H^s(\Omega),$  we have $$ c_1  \|\Gamma_s u\| _{s-1/2 ,\partial \Omega}  \leq \|u\|_{s,\Omega} \leq c_2 \|\Gamma_s u\|_{s-1/2 ,\partial \Omega}, $$
where $\|.\|_{s,\Omega}$ and $\|.\|_{s,\partial \Omega}$ denote the norms  on $H^s(\Omega)$ and $H^s(\partial \Omega)$ respectively.
We may therefore equip $\mathcal H^s(\Omega)$ for $1 \leq s < 3/2$  with the  norm
  $$ \|u\|_{s} =  \|\Gamma_s u\|_{s-1/2,\partial \Omega},$$
 associated with the inner product
  $$(u,v)_{s} = (\Gamma_s u , \Gamma_s v)_{s-1/2,\partial \Omega}.$$
\begin{lemma}
\label{prop:perm} 
Let $\Gamma \in {\mathcal B}( H_{\partial}^1(\Omega),L^2(\partial \Omega))$ 
be the trace operator and let $\Lambda \in {\mathcal C}(L^2(\partial \Omega), 
H_{\partial}^1(\Omega))$ be its Moore--Penrose inverse. Then, for  real $s\geq 0$, 
the following equality holds:
$$
T_{\Lambda^*} (I+\Lambda \Lambda^*)^{-s} 
= (I+\Lambda^* \Lambda)^{-s}  T_{\Lambda^*}, 
$$
where $T_{\Lambda^*}=\Lambda^*(I+\Lambda\Lambda^*)^{-1/2}
+\Gamma(I+\Lambda\Lambda^*)^{-1/2}$. In particular, we have \begin{equation*}
\Lambda^*(I+\Lambda\Lambda^*)^{-1/2}(I+\Lambda \Lambda ^*)^{-s} = (I+\Lambda^* \Lambda )^{-s} \Lambda^*(I+\Lambda\Lambda^*)^{-1/2}.
\end{equation*}
\end{lemma}

\begin{proof}
Since
\begin{equation*}
\begin{split}
\mathcal N(T_{\Lambda^*} (I+\Lambda \Lambda^*)^{-s}) 
&= \mathcal N(T_{\Lambda^*})\\
&= \mathcal N \left(\Lambda^*(I+\Lambda \Lambda^*)^{-1/2} \right)\\
&= \mathcal N( \Lambda^*)\\
&= H_0^1(\Omega)\\
&=\mathcal N \left((I+\Lambda^* \Lambda)^{-s} T_{\Lambda ^*} \right),
\end{split}
\end{equation*}
 we have  
$$
T_{\Lambda^*} (I+\Lambda \Lambda^*)^{-s} v  = 0
= (I+\Lambda^* \Lambda)^{-s}  T_{\Lambda^*} v
$$
for all $v\in \mathcal N(\Lambda^*)$. Now let us consider the operator 
$\Gamma^* \Gamma:  H_{\partial}^1(\Omega) \longrightarrow  H_{\partial}^1(\Omega)$,
where $\Gamma^*$ is the adjoint of the trace operator $\Gamma,$ it is self-adjoint and compact according to the compactness of $\Gamma.$
 Then there exists a sequence of pairs 
$((s_k,v_k))_{k\geq 1}$ associated to $\Gamma^* \Gamma$ such that 
$
\Gamma^*\Gamma v_k = s_k^2 v_k.
$  
To prove the equality on the orthogonal complement $\mathcal N(\Lambda^*)^{\perp}$, 
we shall show that   
$$ 
T_{\Lambda ^*}( I+\Lambda \Lambda^*)^{-s} v_k 
=  (I+\Lambda^* \Lambda)^{-s} T_{\Lambda ^*} v_k.
$$
To this end, let $\Gamma v_k = s_k z_k$. It follows that 
$\Gamma^* z_k = s_k v_k$ and, from Theorem~\ref{thm:3.8}, 
$\Gamma v_k = (I+\Lambda ^* \Lambda)^{-1/2} T_{\Lambda ^*} v_k$.
On the other hand,
$$ 
T_{\Lambda^*} v_k = \Lambda ^* (I+\Lambda \Lambda^*)^{-1/2} v_k 
+ \Gamma (I+\Lambda \Lambda^*)^{-1/2} v_k. 
$$
By putting $(I+\Lambda \Lambda ^*)^{-1}v_k = w_k$, we have 
$v_k = w_k +\Lambda \Lambda ^* w_k$ and 
$$
\Gamma ^* \Gamma v_k = s_k^2 v_k 
= \Gamma ^* \Gamma w_k +w_k,
$$
which implies that 
\begin{equation}
\label{eq:poft1:n}
\begin{split}
(I+\Gamma^* \Gamma)^{-1} \Gamma^* \Gamma v_k 
&= s_k^2 (I+\Gamma^* \Gamma)^{-1} v_k \\
&= w_k \\
&= \Gamma^* \Gamma (I+\Gamma^* \Gamma)^{-1} v_k.
\end{split}
\end{equation}
It follows that 
$$
(I+\Gamma^* \Gamma)^{-1} \Gamma^* \Gamma v_k 
= \Gamma^* \Lambda ^* (I+\Lambda \Lambda^*)^{-1} v_k.
$$
This leads to 
\begin{equation*}
\begin{split}
(I+ \Lambda \Lambda ^*)^{-1} v_k 
&= s_k^2 (I+ \Gamma^* \Gamma )^{-1} v_k \\
&= s_k^2(v_k -(I+\Lambda \Lambda^*)^{-1} v_k),
\end{split}
\end{equation*}
so that 
$$ 
(1+s_k^2) (I+\Lambda \Lambda^*)^{-1} v_k 
= s_k^2 v_k. 
$$
Thus,
\begin{equation*}
\begin{split}
(I+ \Lambda \Lambda^*)^{-1} v_k 
&= \frac{s_k^2}{1+s_k^2}
v_k \\
&= (1+s_k^2)^{-1} s_k^2 v_k,
\end{split}
\end{equation*}
which implies that  
$$ 
\left(I+ \Lambda \Lambda^*\right)^{-s} v_k 
= \left(s_k^2\left(1+s_k^2\right)^{-1}\right)^{s} v_k.
$$
In particular,
$$ 
(I+ \Lambda \Lambda^*)^{-1/2} v_k 
= s_k(1+s_k^2)^{-1/2} v_k.
$$
Consequently,
\begin{equation*}
\begin{split}
T_{\Lambda^*} v_k 
&=  \Lambda^*(I+\Lambda \Lambda^*)^{-1/2} v_k 
+ \Gamma(I+\Lambda \Lambda ^*)^{-1/2} v_k \\
&= \Lambda ^* \left(  \frac{s_k}{\sqrt{1+s_k^2}} v_k \right)
+ \Gamma \left( \frac{s_k}{\sqrt{1+s_k^2}} v_k \right) \\
&= \frac{1}{\sqrt{1+s_k^2}} z_k + \frac{s_k^2}{\sqrt{1+s_k^2}} z_k \\
&= \sqrt{1+s_k^2} z_k.
\end{split}
\end{equation*}
Therefore,
\begin{equation*}
\begin{split}
T_{\Lambda^*}(I+\Lambda \Lambda^*)^{-s} v_k 
&= \left( \frac{s_k^2} {1+s_k^2}\right)^s 
\left( \frac{1+s_k^2}{\sqrt{1+s_k^2}}\right) z_k \\
&= \left(\frac{s_k^2} {1+s_k^2} \right)^s \sqrt{1+s_k^2} z_k.
\end{split}
\end{equation*}
% aqui
On the other hand, 
$\Gamma \Gamma^* z_k = \Gamma s_k v_k = s_k^2 z_k$.
By putting  
$(I+\Lambda^*  \Lambda)^{-1} z_k = e_k$, we have 
$$ 
z_k = e_k+\Lambda ^* \Lambda e_k,
$$
which implies that 
$$ 
\Gamma \Gamma^* z_k = s_k^2 z_k = \Gamma \Gamma^* e_k +e_k 
$$
and 
\begin{equation*}
\begin{split}
(I+\Gamma \Gamma^*)^{-1} \Gamma \Gamma^* z_k 
&= s_k^2 (I+\Gamma \Gamma^*)^{-1} z_k \\
&=e_k\\
&=\Gamma \Gamma^* (I+\Gamma \Gamma^*)^{-1} z_k. 
\end{split}
\end{equation*}
 It follows that 
$$
(I+\Gamma \Gamma^*)^{-1} \Gamma \Gamma^* z_k = \Gamma \Lambda (I+\Lambda^*  \Lambda)^{-1} z_k.
$$ 
This leads to 
$$
(I+\Lambda^*  \Lambda)^{-1} z_k = s_k ^2(I+\Gamma \Gamma ^*)^{-1} z_k 
= s_k^2 \left(z_k - (I+\Lambda^*  \Lambda)^{-1} z_k \right),
$$
so that 
$$
(1+s_k^2) (I+\Lambda^*  \Lambda)^{-1}z_k = s_k^2 z_k 
$$
and
$$    
(I+\Lambda^*  \Lambda)^{-1} z_k = s_k^2 (1+s_k^2) ^{-1} z_k. 
$$
Consequently, 
$$
\left(I+\Lambda^*  \Lambda\right)^{-s} z_k 
= \left(s_k^2 \left(1+s_k^2\right)^{-1}\right)^s z_k, 
$$
which implies that 
\begin{equation*}
\begin{split}
\left( I+ \Lambda^* \Lambda \right) ^{-s} T_{\Lambda^*} v_k 
&=  \left( I+\Lambda ^* \Lambda \right) ^{-s} \sqrt{1+s_k^2} z_k \\
&= \sqrt{1+s_k^2} \left( \frac{s_k^2}{1+s_k^2} \right) ^s z_k. 
\end{split}
\end{equation*}
Hence, one has 
$ \left( I+ \Lambda^* \Lambda \right) ^{-s} T_{\Lambda^*} v_k  
=  T_{\Lambda^*}\left( I+ \Lambda \Lambda^* \right) ^{-s} v_k$
for all $k\geq 1.$ In particular, since $T_{\Lambda^*}=\Lambda^*(I+\Lambda\Lambda^*)^{-1/2}
+\Gamma(I+\Lambda\Lambda^*)^{-1/2},$ we have
\begin{equation*}
\Lambda^*(I+\Lambda\Lambda^*)^{-1/2}(I+\Lambda \Lambda ^*)^{-s} = (I+\Lambda^* \Lambda )^{-s} \Lambda^*(I+\Lambda\Lambda^*)^{-1/2}.
\end{equation*}
\end{proof}
In a similar way, one can prove that:
\begin{lemma}
\label{co:perm}
Let $\Gamma \in {\mathcal B}( H_{\partial}^1(\Omega),L^2(\partial \Omega))$ 
be the trace operator and let $\Lambda \in {\mathcal C}(L^2(\partial \Omega), 
H_{\partial}^1(\Omega))$ be its Moore--Penrose inverse. Then, for real $s\geq 0$, 
the following equality holds:
\begin{equation*}
\Lambda(I+\Lambda ^*\Lambda)^{-1/2}(I+\Lambda^* \Lambda )^{-s} = (I+\Lambda \Lambda ^* )^{-s} \Lambda(I+\Lambda ^*\Lambda)^{-1/2}.
\end{equation*}
\end{lemma}

Consider now for $ v\in \mathcal H_{s}(\Omega)$ the norm $$ \|v\|_{*,s} = \| (I+\Lambda \Lambda ^*)^{s-1} v \| _{\partial,\Omega}, $$
associated with the inner product
 $$ (u,v)_{*,s} = \big ( (I+\Lambda \Lambda ^*)^{s-1} u , (I+\Lambda \Lambda ^*)^{s-1} v \big )_{\partial ,\Omega}, $$
where  $\|.\|_{\partial,\Omega}$ and $(.,.)_{\partial,\Omega}$ denote the norm and the inner product on $H_{\partial}^1(\Omega)$ respectively.  Now we are able to prove Theorem \ref{thm:first main result}.
  \subsubsection{\textbf{Proof of Theorem \ref{thm:first main result}}}
 Assume that $1 \leq  s<3/2.$ Consider $v\in \mathcal H^s(\Omega),$ and let 
   $\Gamma_s v \in H^{s-1/2}(\partial \Omega)$ its trace, and let $\widehat{\Gamma}_s v$  be the embedding of $\Gamma_s v$ in $L^2(\partial \Omega).$ It follows that $\tilde{v}= \Lambda \widehat{\Gamma}_s v$ is the embedding of $v$ in $\mathcal H^1(\Omega).$  Now, since $\widehat{\Gamma}_s v \in H^{s-1/2}(\partial \Omega)$ and that $H^{s-1/2}(\partial \Omega) =\mathcal H^{s-1/2}(\partial \Omega)$ with equivalence of norms (See Theorem \ref{thm:trace}), it follows that there exists $g \in L^2(\partial \Omega)$ such that
   \begin{equation}
   \label{emb}
    \widehat{\Gamma}_s v = \left ( I+\Lambda^* \Lambda \right)^{-(s-1/2)}g,
    \end{equation}
which implies  that
  $$\tilde{v}=\Lambda \widehat{\Gamma}_s v = \Lambda (I +\Lambda^* \Lambda)^{-(s-1/2)} g = \Lambda (I+\Lambda ^* \Lambda)^{-1/2} (I+\Lambda^*  \Lambda)^{-(s-1)}   g ,$$
and from Lemma \ref{co:perm}, we obtain  that
   \begin{equation}
   \label{tilde}
    \tilde{v}  = (I+\Lambda  \Lambda^*)^{-(s-1)} \Lambda (I+\Lambda ^* \Lambda)^{-1/2} g .
    \end{equation}
 Now by putting
  $ w = \Lambda (I+\Lambda ^* \Lambda)^{-1/2}g$ which is   in $\mathcal H^1(\Omega),$ it follows that  $$ \tilde{v}= (I+\Lambda  \Lambda^*)^{-(s-1)} w \in  \mathcal H_{s}(\Omega),$$ we therefore get the first inclusion $\mathcal H^s(\Omega) \subset \mathcal H_s(\Omega).$
  In order to prove the second one,  we consider $v\in \mathcal H_{s}(\Omega).$ There exists then $w\in \mathcal H^1(\Omega)$ such that
  $$v = (I+\Lambda \Lambda^*) ^{-(s-1)} w .$$

  Since $\Gamma = (I+\Lambda ^* \Lambda)^{-1/2} T_{\Lambda ^*}$ (see Theorem~\ref{thm:3.8}), we have 
 
\begin{equation*}
 \Gamma v = (I+\Lambda ^* \Lambda)^{-1/2}  T_{\Lambda^*} (I+\Lambda  \Lambda^* )^{-(s-1)} w 
\end{equation*}
and from Lemma \ref{prop:perm}, we obtain that 
 \begin{equation*}
\begin{split}
  \Gamma v 
   &= (I+\Lambda ^* \Lambda)^{-1/2}  (I+\Lambda^*  \Lambda )^{-(s-1)} T_{\Lambda^*} w \\
   &= (I+\Lambda ^* \Lambda )^{-1/2-s +1} T_{\Lambda^*} w\\
&=  (I+\Lambda ^* \Lambda )^{-(s-1/2)} T_{\Lambda^*} w.
\end{split}
\end{equation*}
 Consequently $$ \Gamma v \in \mathcal H^{s-1/2}(\partial \Omega) = H^{s-1/2}(\partial \Omega).$$
 The uniqueness of the solution of the Dirichlet problem for the Laplace's equation in $H^1(\Omega)$ implies that $v\in \mathcal H^s(\Omega).$ Hence, the algebraic equality $\mathcal H^s(\Omega)= \mathcal H_s(\Omega)$ holds. To prove the equivalence of norms, let us consider $v\in \mathcal H^s(\Omega).$ As in \eqref{emb},  there exists $g \in L^2(\partial \Omega)$ such that
  $$ \widehat{\Gamma}_s v = \left ( I+\Lambda^* \Lambda \right)^{-(s-1/2)}g,$$
which implies that $$\mathcal R(\widehat{\Gamma}_s) \subset \mathcal R(\left( I+\Lambda^* \Lambda \right)^{-(s-1/2)}) ,$$ and according to Theorem \ref{thm:DT}, there exists a bounded operator $T_s$ from $\mathcal H^s(\Omega)$ to $ L^2(\partial \Omega)$ such that
 $$ \widehat{\Gamma}_s v =\left( I+\Lambda^* \Lambda \right)^{-(s-1/2)} T_s v ,$$
 which implies that 
 \begin{equation*}
  \tilde{v} =\Lambda \widehat{\Gamma}_s v = \Lambda \left( I+\Lambda^* \Lambda \right)^{-(s-1/2)} T_s v,
 \end{equation*}
 and  according to Lemma \ref{co:perm}, we obtain that
  \begin{equation}
  \label{tilde1}
   \tilde{v}  =  \left( I+\Lambda \Lambda^* \right)^{-(s-1)} \Lambda \left( I+\Lambda^* \Lambda \right)^{-1/2} T_s v.
  \end{equation}
On the other hand, it follows  from \eqref{tilde} that
 \begin{equation*}
 \begin{split}
 \|\tilde{v}\|_{*,s} &= \| \left( I+\Lambda \Lambda^* \right)^{s-1} \tilde{v} \|_{\partial, \Omega}  \\
 &=  \| \Lambda \left( I+\Lambda^* \Lambda \right)^{-1/2} g\| _{\partial, \Omega}. 
 \end{split}
  \end{equation*}
Hence from  \eqref{tilde1}, we have
  \begin{equation*} 
\|\tilde{v}\|_{*,s} =\| \Lambda \left( I+\Lambda^* \Lambda \right)^{-1/2} g\| _{\partial, \Omega} = \| \Lambda \left( I+\Lambda^* \Lambda \right)^{-1/2} T_s v\| _{\partial, \Omega}.
 \end{equation*}
Using Lemma \ref{lem:LM}, one has
\begin{equation*}
\|T_s v\|_{0, \partial \Omega}^2= \|\Lambda(I+\Lambda^*\Lambda)^{-1/2} 
T_s v\|_{\partial,\Omega}^2 + \|(I+\Lambda^*\Lambda)^{-1/2} T_s v\|_{0, \partial \Omega }^2,
\end{equation*}
which implies that
\begin{equation*}
\|\Lambda(I+\Lambda^*\Lambda)^{-1/2} 
T_s v\|_{\partial,\Omega} \leq \|T_s v\|_{0, \partial \Omega}.
\end{equation*}
The boundedness of $T_s$ implies then that there exists $c>0$ such that
\begin{equation*}
 \|\tilde{v}\|_{*,s}  \leq c \ \|v\|_{s}.
\end{equation*}
Therefore according to the Banach's bounded inverse theorem (see \cite{Co1}), the norms $\|\cdot\|_s$ and $\|\cdot\|_{*,s}$ are equivalent.

  \subsection{\textbf{The case $0\leq s\leq 1$}}

Consider  the embedding operator $E$ of $H_{\partial}^1(\Omega)$ in $L^2( \Omega)$ and its adjoint operator $E^*.$ 
For $f\in L^2(\Omega),$  $E^*$  is the solution operator of the following Robin problem for the Poisson equation

  \begin{equation}
  \label{eq:7.1}
          \begin{cases}
            -\Delta u=f & \text{  }  (\Omega) \\
            \partial_{\nu}u+\Gamma u =0 & \text{}  (\partial \Omega),
          \end{cases}
          \end{equation}
and denote by $E_0^*$  the solution operator of the following Dirichlet problem for the Poisson equation
 \begin{equation}
 \label{eq:7.2}
              \begin{cases}
                -\Delta u^0=f & \text{  }  (\Omega) \\
                \Gamma u^0=0 & \text{}  (\partial \Omega).
              \end{cases}
              \end{equation}

By setting $E_1^*= E^*-E_0^*$ and $u^1 =E_1^*f,$ it follows that $u^1$ is the solution of  the following Dirichlet problem for the Laplace's equation
 \begin{equation}
 \label{eq: DPL}
              \begin{cases}
                \Delta u^1=0 & \text{  }  (\Omega) \\
               \Gamma u^1=\Gamma u & \text{}  (\partial \Omega),
              \end{cases}
              \end{equation}
   where  $u$ is the solution of the problem \eqref{eq:7.1}.
    A consequence of Theorem~\ref{thm:MC} 
   is the following  Rellich--Ne\v{c}as lemma \cite{N}:
     \begin{lemma}[ See Chapter 5 of \cite{N}]
     \label{lem: RN}
   Let $f\in L^2(\Omega)$ and $u^0=E_0^*f$ be the solution of the Dirichlet problem 
   for the Poisson equation \eqref{eq:7.2}. Then, $\partial _{\nu} u^0 \in L^2(\partial \Omega)$. 
   Moreover, there exists a constant $c_{\Omega}>0$, depending on the geometry of $\Omega$, such that
   $$  
   \|\partial_{\nu} u^0 \|_{0,\partial \Omega} \leq  c_{\Omega} \ \|f\|_{0,\Omega}.
   $$
   \end{lemma}
\begin{remark} We can also proceed like in \cite{ch} to prove Lemma \ref{lem: RN}.
\end{remark}
Based on   Rellich--Ne\v{c}as lemma \ref{lem: RN}, it makes sens to define for $g\in L^2(\partial \Omega)$ a very weak solution for the following Dirichlet problem for the Laplace's equation  
     \begin{equation}
     \label{eq:7.4}
             \begin{cases}
               \Delta v=0 & \text{  }  (\Omega) \\
               v =g & \text{}  (\partial \Omega),
             \end{cases}
             \end{equation}
             as follows. We say that $v$ is a very weak solution of the problem \eqref{eq:7.4} if for all $u$ in $H_{\Delta}^1(\Omega) \cap H_0^1(\Omega),$ we have 
                         \begin{equation}
                         \label{very weak}
                         - \int_{\Omega} v \ \Delta u \ dx + \int_{\partial \Omega} g \  \partial _{\nu} u \ d\sigma  =0. 
                         \end{equation}
               
   Denote by   $$
            \begin{array}{ccccc}
            K &: & L^2(\partial \Omega) & \to & L^2(\Omega) \\
             & & g & \mapsto & v=Kg \\
            \end{array}$$ its solution operator.  The adjoint operator $K^*,$  takes each $f \in L^2( \Omega)$ to $-\partial_{\nu} u^0$ into $L^2(\partial \Omega),$ where $u^0$ is the solution of the problem \eqref{eq:7.2}. Indeed, 
            from \eqref{very weak}, we have for all $u$ in $H_{\Delta}^1(\Omega) \cap H_0^1(\Omega),$
        
    $$ - \int_{\Omega} v \ \Delta u \ dx = - \int_{\partial \Omega} g \  \partial _{\nu} u \ d\sigma $$
    and by putting $f= -\Delta u,$ we have
       \begin{equation*}
               \begin{split}
    - \int_{\Omega} v \ \Delta u \ dx    &=   \int_{\Omega} f \ Kg \ dx \\ &=\int_{\partial \Omega} K^* f \ g \ d\sigma\\ &= - \int_{\partial \Omega} \  \partial _{\nu} u \  g  \ d\sigma . 
      \end{split}
      \end{equation*}
       
            Now we will prove Theorem \ref{thm: regularity for K}.
             
  \subsubsection{\textbf{Proof of Theorem \ref{thm: regularity for K}}}
  Let  $f\in L^2(\Omega),$   $ E_1^*f= E^*f- E_0^*f $ and  $u=E^*f$ the solution of the problem \eqref{eq:7.1}.
  Since $\partial_{\nu}u = - \Gamma u \in L^2(\partial \Omega)$ and $\Delta u \in L^2(\Omega),$ then it follows from Theorem~\ref{thm:MC} that $\Gamma u \in H^1(\partial \Omega).$ As $\Gamma u^1 = \Gamma u,$ we have $ \Gamma u^1 \in H^1(\partial \Omega),$  which implies according to  Theorem~\ref{thm:MC} that $\partial_{\nu} u^1 \in L^2(\partial \Omega).$
 If we set $y= \partial_{\nu}u^1 +\Gamma u^1 \in L^2(\partial \Omega),$ then $u^1$ is the solution of the following problem:
    \begin{equation*}
            \begin{cases}
              \Delta u^1=0 & \text{  }  (\Omega) \\
              \partial_{\nu}u^1+\Gamma u^1 =   y & \text{}  (\partial \Omega),
            \end{cases}
            \end{equation*}

  which means that $$u^1 = \Gamma ^* y = E_1^* f,$$
 so$$\mathcal R(E_1^*) \subset \mathcal R(\Gamma ^*). $$ It follows according to Douglas theorem (Theorem~\ref{thm:DT}) that there exists an operator $T\in \mathcal B(L^2(\Omega), L^2(\partial \Omega))$ with $\mathcal N(T) =\mathcal N(E_1^*)$ such that $  E_1^* = \Gamma ^* T.$
  Applying Green's formula (Corollary~\ref{prop:GF} ), we have for all $v\in \mathcal  H^1(\Omega)$  
  \begin{equation*}
  \begin{split}
  (u^1,v)_{\partial ,\Omega}&= \int_{\Omega} \nabla u^1 \nabla v \ dx + \int_{\partial \Omega} \Gamma u^1 \Gamma v \ d\sigma \\
  &= - \int _{\Omega} \Delta u^1 Ev \ dx + \int_{\partial \Omega}\partial _{\nu} u^1 \Gamma v \ d\sigma + \int_{\partial \Omega} \Gamma u^1 \Gamma v \ d\sigma \\
  &= \int_{\partial \Omega} (\partial_{\nu}u^1 + \Gamma u^1) \ \Gamma v \ d\sigma ,
  \end{split}
  \end{equation*}
  so that if we set $g= \partial_{\nu}u^1 + \Gamma u^1,$  we have
   $$(u^1,v)_{\partial, \Omega}=(g,\Gamma v)_{0,\partial \Omega}= (\Gamma^*g,v)_{\partial,\Omega}.  $$
  Since $v$ is arbitrary in $\mathcal H^1(\Omega),$ we deduce that $u^1 = \Gamma^* g $  and this implies that for $f\in L^2(\Omega),$ 
  $$Tf=  \partial_{\nu}u^1 + \Gamma u^1.$$
  On the other hand, $u^1,u^0 \in H_{\Delta}^1(\Omega)$ implies that 
  $\widehat{\partial}_{\nu} u^1, \widehat{\partial}_{\nu} u^0 \in H^{-1/2}(\partial \Omega),$
  so $$\widehat{\partial}_{\nu}u = \widehat{\partial}_{\nu} u^0 +\widehat{\partial}_{\nu} u^1 \in  H^{-1/2}(\partial \Omega).$$ Now, since $\partial_{\nu} u,\partial_{\nu} u^1 \in L^2(\partial \Omega), $ it follows that $\partial_{\nu} u^0 \in L^2(\partial \Omega)$
  and that $$ - \partial_{\nu} u^0 = \partial_{\nu} u^1  +\Gamma u^1 = Tf.$$
  On the other hand,   $K^*$ maps each $f\in L^2(\Omega)$ to $-\partial_{\nu} u^0 $ onto $L^2(\partial \Omega),$  where $u^0$ is the solution of Dirichlet problem for the Poisson equation \eqref{eq:7.2}.
  This leads to
  $$ - \partial_{\nu} u^0 = K^*f,$$
  therefore  $$T= K^*. $$
  The proof is complete.
\subsubsection{\textbf{An orthonormal basis for the $L^2-$ Bergman space}}
The operator $E_1$ being bounded, one considers its Moore-Penrose inverse which we denote by $F_1$ such that $$\mathcal D(F_1) = \mathcal R(E_1) \oplus \mathcal N(E_1^*),$$
and 
  $$\mathcal R(F_1)=\mathcal H^1(\Omega) \ \ \mbox{and} \ \ \mathcal N(F_1)= \mathcal N(E_1^*).$$
From Theorem ~\ref{thm:3.5} 
 the operator $ F_1^* (I+F_1 F_1^*)^{-1/2}$ from $ H^1(\Omega)$ into $L^2(\Omega)$  is bounded and has a closed range such that 
$$\mathcal R  (F_1^* (I+F_1 F_1^*)^{-1/2}) = \mathcal R (F_1^*)= \mathcal H(\Omega)$$
and
$$\mathcal N( F_1^* (I+F_1 F_1^*)^{-1/2}) = \mathcal N (F_1^*) = H_0^1(\Omega) .$$ Moreover,
 from Corollary~\ref{cor:3.6}, we have:
 \begin{lemma}
 \label{lem:7.6}
  The operator
    $ F_1^* (I+F_1 F_1^*)^{-1/2} \in \mathcal B(\mathcal H^1(\Omega), \mathcal H(\Omega))$ is an isomorphism.
 \end{lemma} 
  \label{lem:7.7}
Let us now set $$\Gamma_0^* = F_1^* (I+F_1 F_1^*)^{-1/2} \Gamma ^* ,$$
where $\Gamma^*$ is the adjoint operator of $\Gamma \in \mathcal B (H_{\partial}^1(\Omega), L^2(\partial \Omega)).$
\begin{lemma}
The operator $\Gamma_0^*$ defined above  is compact and injective.
\end{lemma}
\begin{proof}
Knowing that $\mathcal R  (F_1^* (I+F_1 F_1^*)^{-1/2}) = \mathcal H(\Omega),$
we have $ \mathcal R(\Gamma_0^*) \subset \mathcal H(\Omega). $ Moreover, the operator $\Gamma^*$ is compact according to Lemma \ref{lem:TA}.  Therefore, the boundedness of $F_1^* (I+F_1 F_1^*)^{-1/2} $ implies  that  $\Gamma_0^*$ is compact. The injectivity of $\Gamma_0^*$ holds for the reason that $\Gamma^*$ is injective (see Lemma \ref{lem:TA}) and that $\mathcal R(\Gamma^*)\subset \mathcal N(F_1^* (I+F_1 F_1^*)^{-1/2})^{\perp}= \mathcal H^1(\Omega),$ where $\mathcal N(F_1^* (I+F_1 F_1^*)^{-1/2})^{\perp}$ is the orthogonal complement of $\mathcal N(F_1^* (I+F_1 F_1^*)^{-1/2}).$
\end{proof}
Theorem \ref{thm:comp}  provides an orthonormal basis for the Bergman space $\mathcal H(\Omega).$ The proof is as follows.

  \subsubsection{\textbf{Proof of Theorem \ref{thm:comp}}}  
   By definition we have  $$\Gamma_0^* = F_1^* (I+F_1 F_1^*)^{-1/2} \Gamma ^* .$$   Composing  by $K^*,$ we obtain that
          \begin{equation*}
          \begin{split}
          \Gamma_0^*  K^* =  F_1^* (I+F_1 F_1^*)^{-1/2} \Gamma ^* K^* , 
          \end{split}
          \end{equation*}
        and in view of Theorem~\ref{thm: regularity for K}, we have $E_1^* = \Gamma^* K^*,$ which implies that
         \begin{equation*}
           \Gamma_0^*  K^* = F_1^* (I+F_1 F_1^*)^{-1/2} E_1^*  .
           \end{equation*}
           On the other hand, since $(I+F_1 ^* F_1)^{-1/2} F_1^* \subset  F_1^* (I+F_1 F_1^*)^{-1/2}$ and that  $\mathcal R(E_1^*)\subset \mathcal D(F_1^*),$ one has
          $$
           \Gamma_0^*  K^*= (I+F_1 ^* F_1)^{-1/2} F_1^* E_1^* , $$
           and viewing $F_1^* E_1^* =P_{\mathcal H(\Omega)},$ it follows that
           $$
         \Gamma_0^*  K^*=(I+F_1 F_1^*)^{-1/2} P_{\mathcal H(\Omega)}.$$
          On the other hand, viewing  $\Gamma_0^*$ is compact,  $K\in \mathcal B(L^2(\partial \Omega),L^2(\Omega))$   and that 
                     $ (I+F_1^* F_1) ^{-1/2}  P_{\mathcal H(\Omega)} $ is  a positive self-adjoint, it follows that $\Gamma_0^* K^*$  is a compact self-adjoint operator. Therefore, according to the spectral theorem for compact self-adjoint operators we have the existence of an orthonormal basis  $(\phi_n)_{n\in\mathbb N}$ for $\mathcal H(\Omega)$   and a real sequence $(\kappa_n)_n$ such that $\displaystyle{\lim_{n \mapsto {+ \infty}}} \kappa_n =0$  and for all $n\in \mathbb N,$ $ \Gamma _0^* K^* \phi _n =  \kappa _n^2 \phi_n .$ The proof is complete.

 \subsubsection{\textbf{Proof of Theorem \ref{thm: main result}}}
For $s=0,$ the equality 
$\mathcal H(\Omega)= \mathcal X^0(\Omega)$ holds 
by definition. For $s=1,$ we have from Theorem~\ref{thm:3.8} the following decomposition $$E_1=(I+F_1^* F_1)^{-1/2}T_{F_1^*}$$
where $$T_{F_1^*}  = F_1^*(I+F_1F_1^*)^{-1/2} + E_1(I+F_1F_1^*)^{-1/2} .$$ 
 By definition we have $\mathcal R (E_1 )= \mathcal H^1(\Omega).$ On the other hand, since $T_{F_1^*}$ is an isomorphism of $\mathcal H^1(\Omega)$ and  $\mathcal H(\Omega)$ according to Corollary \ref {cor:2.8},   it follows that 
$$\mathcal R (E_1 )=\mathcal R \left ((I+F_1^* F_1)^{-1/2} P_{\mathcal H(\Omega)} \right )= \mathcal X^1(\Omega),  $$
which implies that $\mathcal H^1(\Omega)= \mathcal X^1(\Omega).$ To prove the equivalence of norms,
 consider for $u\in \mathcal X^1(\Omega)$ the  norm $$\|u\|_{\mathcal X^1(\Omega) }= \|(I+F_1^* F_1)^{1/2}   u\| _{0,\Omega} $$ and for  $v \in \mathcal H^1(\Omega),$
  the norm
  $$ \|v\|_{\mathcal H^1(\Omega)}= \|v \|_{\partial, \Omega}.$$
For $v\in \mathcal H^1(\Omega),$ we have $E_1 v \in \mathcal X^1(\Omega)$ and
$$ \|(I+F_1^* F_1)^{1/2} E_1 v\| _{0,\Omega} = \|T_{F_1^*} v\|_{0,\Omega}. $$
Viewing $T_{F_1^*} $ is an isomorphism of $\mathcal H^1(\Omega)$ and $\mathcal H(\Omega)$ according to Corollary~\ref{cor:2.8}, there exist then two positive  constants $c_1^{\prime}$ and $c_2^{\prime}$ not depending on $v$ such that
$$ c_1^{\prime}  \|v\|_{\partial,\Omega} \leq  \|(I+F_1^* F_1)^{1/2} E_1  v\| _{0,\Omega}  \leq c_2^{\prime} \|v\|_{\partial,\Omega} .$$ Therefore,  the norms  $\|.\|_{\mathcal X^1(\Omega)}$ and  $\|.\|_{\partial ,\Omega}$ are equivalent,  which means that  $\mathcal X^1(\Omega)=\mathcal H^1(\Omega)$ with equivalence of norms.
 Moreover, $\mathcal X^s(\Omega)$ form an interpolating family  according to  classical results on the theory of positive self-adjoint operators (see \cite{Ad}, \cite{Mc} and \cite{T}). We therefore deduce that $ \mathcal H^s(\Omega) =  \mathcal X^s(\Omega)$ with equivalence of norms. Let us now consider $v\in \mathcal X^s(\Omega).$  There exists then $\phi \in \mathcal H(\Omega)$ such that
    $$ v = (I+F_1^*F_1)^{-s/2}\phi ,$$
    which implies  that for all $n \geq 1,$ we have
    $$(v,\phi_n)_{0,\Omega} = ((I+F_1^* F_1)^{-s/2}\phi , \phi_n)_{0,\Omega} = (\phi,(I+F_1^* F_1)^{-s/2}\phi_n)_{0,\Omega}.$$
  According to the Spectral Theorem  \cite{Co1},  it follows that  for all $n \geq 1,$
     $$(v,\phi_n)_{0,\Omega} = ( \phi , (I+F_1^* F_1)^{-s/2}\phi_n)_{0,\Omega} = (\phi, \kappa_n^{2s} \phi_n)_{0,\Omega},$$
  
  which implies that
     $$ \frac{1}{\kappa_n^{2s}} (v,\phi_n)_{0,\Omega} = (\phi, \phi_n)_{0,\Omega}. $$
We therefore obtain 
     $$ \sum_{n=1}^{\infty} \frac{1}{\kappa_n ^ {4s}} |(v,\phi _n)_{0,\Omega}|^2 = \sum_{n=1}^{\infty} |(\phi,\phi_n)_{0,\Omega}|^2 . $$
     Since for  $\phi \in \mathcal H(\Omega),$
     $ \displaystyle{\sum_{n=1}^{\infty}} |(\phi,\phi_n)_{0,\Omega}|^2 < + \infty, $ it follows  that
       $$  \displaystyle{\sum_{n=1}^{\infty}} \frac{1}{\kappa_n ^ {4s}} |(v,\phi _n)_{0,\Omega}|^2  < + \infty,  \ \ \  \forall v \in \mathcal X^s(\Omega),$$
  hence we conclude that  for $v\in \mathcal H^s(\Omega),$
      $  \displaystyle{\sum_{n}} \frac{1}{\kappa_n ^ {4s}} |(v,\phi _n)_{0,\Omega}|^2 $ converges. The proof is complete.

  \subsection{Boundary formulas for $ \mathcal R(K), \mathcal R(\Gamma_0^*)$  and   $\mathcal H^s(\Omega),$ for $1/2<s<3/2$ } 
    \label{sec:06} According to Weyl's lemma \cite[Theorem 2.3.1]{Mor}, we have $\mathcal H^s(\Omega) \subset \mathcal C^{\infty}(\Omega)$  for all $s \geq 0.$ Also, using the Mean-value theorem \cite[Theorem 2.2.1]{Mor}, the evaluation functionals associated to $\mathcal H^s(\Omega)$ for all $s \geq 0, $ are continuous, which means that  $\mathcal H^s(\Omega)$ are reproducing kernel Hilbert spaces. In this subsection  we will establish  boundary integral formulas for $\mathcal R(K), \mathcal R(\Gamma_0^*)$ and for $\mathcal H^s(\Omega)$ for the range of values $1/2<s<3/2$. 
To this end, let us consider the orthonormal basis $(\phi_n)_{n\geq 1}$  of $\mathcal H(\Omega)$ given in Theorem \ref{thm:comp}. The Bergman kernel  should then, according to \eqref{eq: forme som}, take the form
         \begin{equation*}b(x,y) = \displaystyle \sum_{n=1}^{\infty} \phi_n(x) \phi_n(y), 
        \end{equation*}
          where the convergence here is uniform on all compact in $\Omega \times \Omega.$
        Moreover,
        for all $v\in \mathcal H(\Omega),$ we have  
       $$
      v= \displaystyle {\sum_ {n=1}^{+\infty}} (v,\phi_n)_{0,\Omega} \ \phi_n
      $$  
     which implies that 
           \begin{equation*}
           \begin{split}
             v(x)&= \displaystyle {\sum_ {n=1}^{+\infty}} (v,\phi_n)_{0,\Omega} \ \phi_n(x)\\
             &=  \displaystyle {\sum_ {n=1}^{+\infty}} \left( \int_{\Omega} v(y) \phi_n(y) \ dy \right) \phi_n (x) \\
             &= \int_{\Omega} \left( \displaystyle {\sum_ {n=1}^{+\infty}} \phi_n (x) \phi_n(y) \right) v(y)\  dy \\ 
         & =  \int_{\Omega} b(x,y) \ v(y) \ dy.
         \end{split}
              \end{equation*}
       
    \subsubsection{\textbf{Proof of Theorem \ref{thm:BF}}}
     Let $A$ be a bounded operator from $L^2(\partial \Omega)$ to $  L^2( \Omega) $  such that  $\mathcal R(A) \subset \mathcal H(\Omega)$  and let $v\in \mathcal R(A).$ There exists then $g\in L^2( \partial \Omega)$ such that $v= A g.$
       It follows that
    \begin{equation*}
        \begin{split}
                       Ag(x) &=  \int_{\Omega} b(x,y) Ag(y) \ dy  \\
                       &= \int_{\Omega} b_x(y) Ag(y)
                         \ dy \\
                       &=  \int_{\partial \Omega} A^* b_x(y) g(y)
                                            \ d\sigma(y), 
                                            \end{split} \end{equation*}
 and the proof is complete.

   \begin{remark}
 
 In particular, for $v\in \mathcal H^1(\Omega), $ we have the boundary formula
    \begin{equation}
    \label{eq:5.2}
    E_1v(x)=\int_{\partial \Omega} \left(  K^* b_x \right ) (y) \ \Gamma v(y) d\sigma(y).
    \end{equation}
   \end{remark}

   \subsubsection{\textbf{Proof of Corollary \ref{co:1}}}
    For $ 1/2<s<3/2,$ let us consider $v\in \mathcal H^s(\Omega),$ let  $\Gamma_sv \in H^{s-1/2}(\partial \Omega)$ be its trace, and let  $\widehat{\Gamma}_sv$ be the embedding of $\Gamma_sv$  in $L^2(\partial \Omega).$  It follows that $v$ is the unique very weak solution of 
     \begin{equation*}
                 \begin{cases}
                   \Delta v=0 & \text{  }  (\Omega) \\
           v = \widehat{\Gamma} _s v& \text{}  (\partial \Omega)
                 \end{cases}
                 \end{equation*}
     and that $$v= K \widehat{\Gamma}_s v.$$
       Since $\mathcal R(K \widehat{\Gamma}_s) \subset \mathcal H(\Omega),$ it follows that 
       \begin{equation*}
       \begin{split}
       v(x)&= \int_{ \Omega} b(x,y) K \widehat{\Gamma}_s  v dy \\ 
       &= \int_{\partial \Omega} \left(  K^* b_x \right ) (y) \left( \widehat{\Gamma}_s v \right )(y) \ d\sigma(y).
       \end{split}
       \end{equation*}
     \subsubsection{\textbf{Proof of Corollary \ref{co:2}}}
      Immediate consequence of Theorem \ref{thm:BF} since $\mathcal R(K) \subset \mathcal H(\Omega).$
       \begin{remark}
            The formula \eqref{Lions's formula} corresponds to Lions' formula according to  Englis, Lukkassen, Peetre and Person in \cite{ELPP} where  $\left(  K^* b_x  \right )$ is the Poisson kernel for the Laplace's equation.
         \end{remark}
             \subsubsection{\textbf{Proof of Corollary \ref{co:3}}}
               Since  $\Gamma_0^* = F_1^* (I+F_1 F_1^*)^{-1/2} \Gamma ^* ,$ and $$\mathcal R  (F_1^* (I+F_1 F_1^*)^{-1/2}) = \mathcal R (F_1^*)= \mathcal H(\Omega)$$it follows that $\mathcal R(\Gamma_0^*) \subset \mathcal H(\Omega),$  and according to Theorem \ref{thm:BF}, we deduce that 
                   for all $v\in \mathcal R(\Gamma_0^*),$ there exists $z\in L^2(\partial \Omega)$ such that
                     
                      \begin{equation}
                      v(x)=\int_{\partial \Omega} \Big(  \Gamma_0 b_x \Big ) (y) z(y) d\sigma(y),
                      \end{equation}
                      where  $b(\cdot, \cdot)$ is the Bergman kernel given by \eqref{bergman kernel}.
                     \subsubsection{\textbf{Proof of Theorem \ref{prop:repro}}}
    Let  $(\phi_n)_{n\geq 1}$ be the orthonormal basis for $\mathcal H(\Omega)$ stated in Theorem \ref{thm:comp}.
          Since  $(I+F_1^* F_1)$ is  a positive self-adjoint operator on $\mathcal H(\Omega), $ then for all $v,w \in \mathcal H^s(\Omega)$ when $0 \leq s \leq 1,$ one can naturally state that
         $$(v,w)_{s,\Omega} = \left ( (I+F_1^* F_1)^{s/2}v, (I+F_1^* F_1)^{s/2} w \right) _{0,\Omega}. $$
   In particular, for all $k,l \geq 1$ we have
             \begin{equation*}
             \begin{split}
       (\phi_k,\phi_l)_{s,\Omega}
       & = \left ( (I+F_1^* F_1)^{s/2}\phi_k, (I+F_1^* F_1)^{s/2
       } \phi_l \right) _{0,\Omega}\\
       &= \left ( \frac{1}{\kappa_k^{2s}} \phi_k ,\frac{1}{\kappa_l^{2s}} \phi_l \right)_{0,\Omega},
       \end{split}
               \end{equation*}
               which implies that 
               \begin{equation*}
    (\kappa_k^{2s}\phi_k,\kappa_l^{2s}\phi_l)_{s,\Omega}    = (\phi_k,\phi_l)_{0,\Omega}= \delta _{kl},
        \end{equation*}
        where
        $\delta _{kl}$ is the  Kronecker symbol. Hence  for $0 \leq s \leq 1,$ $\left ((\kappa_n^{2s}\phi_n) \right )_{n\geq 1}$ is an orthonormal basis for  $\mathcal H^s(\Omega).$ Moreover, by putting   $\phi_n^s = \kappa_n^{2s} \phi_n,$
              the associated  reproducing kernel  should take the form:
                $$b^s(x,y) = \displaystyle \sum_{n=1}^{\infty} \phi_n^s(x) \phi_n^s(y). $$ The proof is complete.
  
{\bf Acknowledgments.} This research is part of the second author's Ph.D. dissertation, 
which is carried out at Moulay Ismail University, Mekn\`{e}s. 

% -------------------------------------------

% -------------------------------------------

\end{document}